\documentclass[aps,sd,amsmath,amssymb,reprint,]{revtex4-1}

\draft 

\usepackage{dcolumn}
\usepackage{bm}
\usepackage{graphicx}
\usepackage{caption}
\usepackage[inline]{enumitem}

\usepackage{epstopdf}
\usepackage[english]{babel}
\usepackage{amsthm}
\usepackage[inline]{enumitem} 

\theoremstyle{definition}
\newtheorem{definition}{Definition}

\theoremstyle{remark}

\usepackage{color,soul}

\newcommand{\R}{\mathbb{R}}

\newcommand{\N}{\mathcal{N}}
\newcommand{\A}{\mathcal{A}}
\begin{document}

\title{Weak tracking in nonautonomous chaotic systems} 

\author{Hassan Alkhayuon}
\email[]{hassan.alkhayuon@ucc.ie}
\affiliation{School of Mathematical Sciences, University College Cork, Cork, Ireland. }
\author{Peter Ashwin}
\email[]{p.ashwin@exeter.ac.uk}

\affiliation{Centre for Systems, Dynamics and Control, Department of Mathematics, University of Exeter, Exeter EX4 4QF, UK.}

\date{\today}

\begin{abstract}
Previous studies have shown that rate-induced transitions can occur in pullback attractors of systems subject to ``parameter shifts'' between two asymptotically steady values of a system parameter. For cases where the attractors limit to equilibrium or periodic orbit in past and future limits of such an nonautonomous systems, these can occur as the parameter change passes through a critical rate. Such rate-induced transitions for attractors that limit to chaotic attractors in past or future limits has been less examined. In this paper, we identify a new phenomenon is associated with more complex attractors in the future limit: weak tracking, where a pullback attractor of the system limits to a proper subset of an attractor of the future limit system. We demonstrate weak tracking in a nonautonomous R\"ossler system, and argue there are infinitely many critical rates at each of which the pullback attracting solution of the system tracks an embedded unstable periodic orbit of the future chaotic attractor. We also state some necessary conditions that are needed for weak tracking.

\end{abstract}


\maketitle 

\section{Introduction}

Attractors of nonautonomous (time-varying) dynamical systems that limit to autonomous systems in both past and future time can undergo rate-induced transitions. Many studies and applications of these transitions in such ``parameter shift'' systems assume equilibrium attractors of both limiting systems, see Ref. \cite{Ashwin2011, Wieczorek2011, Perryman2014, Ritchie2016, Ritchie2017, Hoyer-Leitzel2017, Alkhayuon2019,OKeeffe2019}. If there are non-equilibrium attractors for past and future limits these can lead to new phenomena. For example, Kasz\'as {\em et al} \cite{Kaszas2016} study the equation of the forced pendulum with time-dependent amplitude of forcing and show there is an analogy between the behavior of the pullback attractor of the nonautonomous system and the bifurcation diagram of the associated autonomous (or ``frozen system''). The structure of the pullback attractor may be very complex even for parameter values where is no stable chaos. In another paper, Kasz\'as {\em et al} \cite{Kaszas2018} explain the time-dependent topology of the same system and show that it can be described using properties of pullback saddles and their unstable foliations. 

Rate-induced transitions for attractors that limit to various sets in the past are discussed in Alkhayuon and Ashwin \cite{Alkhayuon2018}, where each attractor for the past limit system can be associated with a pullback attractor for the nonautonomous system. For such a system with a branch of exponentially stable attractors, Ref. \cite{Alkhayuon2018} identify a number of rate-induced phenomena: 
\begin{enumerate*}
\item {\em strong tracking}: where a pullback attractor of the system end-point tracks the branch of attractors and limit fully to the attractor of the future limit system;
\item {\em partial tipping}: where certain trajectories of a pullback attractor track the branch but other trajectories tip (i.e. limit to other attractors forward in time);
\item {\em total tipping}: where a whole pullback attractor limits forward in time to an attractor that is not included in the considered branch. 
\end{enumerate*}

An invariant set $M$ is called a {\em minimal invariant set} if it contains no proper invariant subset. Analogously, an attractor $A$ is called a {\em minimal attractor} if it has no proper sub-attractors \cite{Milnor1985}. Chaotic attractors  such as the R\"ossler attractor provide a rich source of attractors that are non-minimally invariant; as they typically contain a dense set of embedded unstable periodic orbits  \cite{Letellier2006}.

Assume we have a parameter shift system that limits forward in time to a system with a non-minimal attractor, or even a minimal attractor that is not minimally invariant. We say there is a {\em weak tracking}, if there is a pullback attractor for the parameter shift system that limits forward in time  to one of the invariant subsets of the future limit attractor. The future limit system needs to have at least one attractor that is non-minimal invariant set, in order for the parameter shift system to exhibit weak tracking. This can be seen on applying \cite[Lemma~II.1]{Alkhayuon2018} which shows that the upper forward limit of a pullback attractor must be invariant with respect to the future limit system. 

In this paper we demonstrate the existence of weak tracking of pullback attractors for parameter shift systems. In Section~\ref{sec:parshift} we define weak tracking for parameter shift systems. In doing so, we use the results on asymptotic behaviour of parameter shift systems from Ref. \cite{Alkhayuon2018}. Section~\ref{sec:bif_Ross} illustrate the phenomena of weak tracking in R{\"o}ssler system \cite{Rossler1976}.  We shift one bifurcation parameter of the system monotonically such that future limit system has always a chaotic R\"ossler attractor, whereas, the past limit system has an attracting equilibrium. We show that there is a dense set of critical rate at each of which the system exhibits weak tracking. Finally, we discuss and conclude in Section~\ref{sec:discussion}. In particular, we note a dimension restriction that must be satisfied for weak tracking to take place - the past limit attractor can have dimension no bigger than the stable manifold of a proper subset of the future limit attractor.

\section {Asymptotic behaviour of parameter shift systems} 
\label{sec:parshift}

A {\em parameter shift system} \citep{Ashwin2017} is a nonautonomous differential equation of the form:
\begin{equation}
\dot{x}=f(x,\Lambda(rt)),
\label{eq:parshift}
\end{equation}
where $x\in\R^n$, $t,r\in\R$, $\Lambda:\R\rightarrow \R$ and $f$ is at least $C^1$ in both arguments. For some $\lambda_-$ and $\lambda_+ \in \R$ with $\lambda_- < \lambda_+$, the parameter shift $\Lambda$ satisfies (i) $\Lambda(\tau) \in (\lambda_-,\lambda_+)$ for all $\tau \in \R$, (ii) $\lim_{\tau \to \pm\infty} \Lambda(\tau)=\lambda_{\pm}$, and (iii)
 $\lim_{\tau \to \pm\infty}  d\Lambda/d\tau = 0$.
 We denote the solution process of \eqref{eq:parshift} with $x(s)=x_0$ by $\Phi(t,s,x_0):=x(t)$. One can understand much of the behaviour of System~\eqref{eq:parshift} by studying the associated autonomous (or frozen) system, which is given by: 
 \begin{equation}
 \dot{x}=f(x,\lambda), 
 \label{eq:parshift_fixed}
\end{equation}  
 where $\lambda$ is time-independent and denote the flow of \eqref{eq:parshift_fixed} by $\phi_{\lambda}(t,x_0):=x(t)$, where $x(0)=x_0$. 

We say a set valued function $\mathcal{M} = \{M_t\}_{t\in\R}$ of $t\in\R$ is a {\em nonautonomous set} for \eqref{eq:parshift} if $M_t$ is nonempty for all $t\in\R$ \cite{Kloeden}. Moreover, $\mathcal{M}$ is called $\Phi$-invariant if $\Phi(t,s,M_s) = M_t$ for all $t,s\in\R$. We say that $\mathcal{M}$ has a property $p$ if and only if $M_t$ has $p$ for all $t\in\R$. 

To study the asymptotic behaviour of nonautonomous sets, note there are several different notions of limit for set valued sequences \cite{set-valued}. More precisely, for a nonautonomous set $\mathcal{M} = \{M_t\}_{t\in\R}$ \cite{Rasmussen2008} one can define the upper forward limit ($M_{+\infty}$) and the upper backward limit ($M_{-\infty}$) of $\mathcal{M}$ as follows:
$$
M_{+\infty} :=\limsup_{t\to\infty} M_t =\bigcap_{\tau >0} \overline{\bigcup _{t\geq \tau} M_t}, 
$$ 
$$
M_{-\infty} :=\limsup_{t\to-\infty} M_t=\bigcap_{\tau >0} \overline{\bigcup _{t\leq \tau} M_t}. 
$$ 
We focus on these upper limits (rather than lower limits) as they capture the asymptotic behaviour in maximal sense. 

Furthermore, we denote the set of asymptotically stable attractors  of \eqref{eq:parshift_fixed} that are parameterised by $\lambda$ by $\mathcal{X}_{\text{as}}$. The set of all exponentially stable attractors $\mathcal{X}_{\text{stab}}$ is a subset of $\mathcal{X}_{\text{as}}$.  We call the boundary of $\mathcal{X}_{\text{stab}}$, $\overline{\mathcal{X}_{\text{stab}}} \setminus \mathcal{X}_{\text{stab}}$, set of bifurcations. One can think of them as subsets of $\R^n \times [\lambda_-,\lambda_+]$.   A continuous set valued function $A(\lambda)\in\mathcal{X}_{\text{as}}$, for all $\lambda \in [\lambda_-,\lambda_+]$, is called a {\em stable path}. If $A(\lambda) \in \mathcal{X}_{\text{stab}}$, for all $\lambda \in [\lambda_-,\lambda_+]$, and its stability is independent of $\lambda$, in the sense that the exponential rate of converging to $A(\lambda)$ is independent of $\lambda$ then we say the path is {\em uniformly stable}, for more details see Ref. \cite{Alkhayuon2018}. A uniformly stable path is called a {\em stable branch} \cite{Alkhayuon2018}. Note that a stable path can include a several stable branches joined at bifurcation points, for an example of a stable path that continues bifurcation points, see Section~\ref{sec:bif_Ross}.

\subsection{Weak tracking of pullback attractors}
\label{subsec:weaktrack}

We define local pullback attractors as in \cite{Alkhayuon2018}. Suppose that $\Phi$ is a process on $\R^n$. A compact and $\Phi$-invariant nonautonomous set $\A$ is called {\em local pullback attractor} if there exists an open set $U$ that contains the upper backward limit of $\A$ and satisfies. 
$$
\lim_{s\rightarrow -\infty} d(\Phi(t,s,U),A_t)=0,
$$
for all $t\in \R$, where $d$ is Hausdorff semi-distance.

Theorem~II.2 shows that for each asymptotically stable attractor $A_-$ for the past limit system there is a local pullback attractor for \eqref{eq:parshift} whose upper backward limit is contained in $A_-$. This pullback attractor depends on the parameter shift $\Lambda$, the rate $r$ as well as the attractor of the past limit system $A_-$. Therefore, we denote the pullback attractor by $\A^{[\Lambda,r,A_-]}$ and it consists of $t$-fibres that are defined as
:
\begin{equation}
A_t^{[\Lambda,r,A_-]} := \bigcap_{\tau>0}\overline{\bigcup_{s\leq\tau}\Phi(t,s,\N_\eta(A_-))}
\label{eq:pbattractor}
\end{equation}

for some $\eta>0$. Note that if $A_-$ is an equilibrium then \cite[Theorem~2.2]{Ashwin2017} shows that the pullback attractor is a single trajectory or so-called {\em pullback attracting solution}.

For a uniformly exponentially stable branch $A(\lambda)$ that contains an attractor of the past limit system $A_-:=A(\lambda_-)$ and for sufficiently small positive $r$, \cite[Theorem~III.1]{Alkhayuon2018} proves that the pullback attractor (\ref{eq:pbattractor}) end-point tracks the branch $A(\lambda)$ . 

This tracking is not guaranteed for large values of $r>0$ or where a stable branch is weakened to a stable path. Rate-induced transitions take place when this tracking breaks.  \citep[Definition~III.1]{Alkhayuon2018} defines different rate-induced transitions between Partial tipping, total tipping and invisible tipping. Here we present a new phenomenon we call {\em weak tracking} that can also lead to transitions.

\begin{definition}
Suppose that $(A(\lambda),\lambda)\subset \mathcal{X}_{\text{as}}$ is a path of asymptotically stable attractors for $\lambda\in[\lambda_-,\lambda_+]$. Define $A_{\pm}:=A(\lambda_{\pm})$ and consider the pullback attractor $\mathcal{A}^{[\Lambda,r,A_-]}$ with past limit $A^{[\Lambda,r,A_-]}_{-\infty}$ that is contained in $A_-$.  We say there is {\em strong tracking} for system \eqref{eq:parshift} from $A_-$ for some $\Lambda$ and $r>0$ if 
$
A_{+\infty}^{[\Lambda,r,A_-]} = A_+.
$
 We say there is {\em weak tracking} if $A_{+\infty}^{[\Lambda,r,A_-]} \subsetneq A_+.$
\end{definition}

Lemma~II.1 from Ref. \cite{Alkhayuon2018} shows that the upper forward limit $A_{+\infty}^{[\Lambda,r,A_-]}$ is invariant with respect to the future limit system. Consequently, in order to exhibit weak tracking the future limit system needs to have an attractor with a proper invariant subset. 

As an example of this behaviour we consider the R\"ossler system \cite{Rossler1976} with embedded unstable periodic orbits that can be the upper forward limit of the pullback attractor for some positive $r$.

 \begin{figure*} 
 \centerline{\includegraphics[width=1.8\columnwidth]{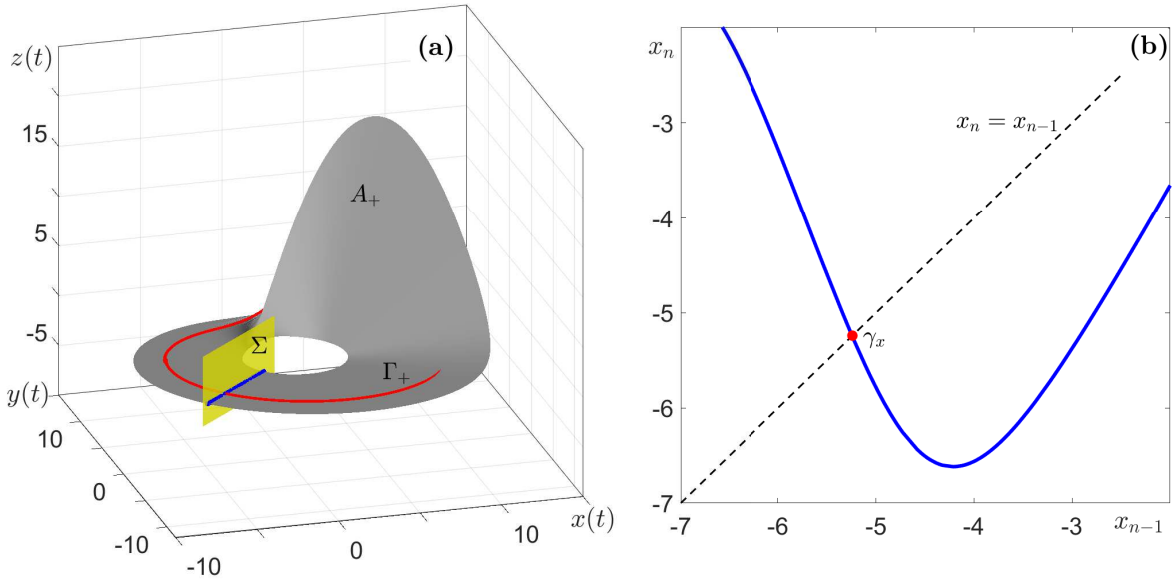}}
   \caption{ In (a), the R\"ossler attractor for parameter values $a=b=0.2$ and $c=5.7$. This also shows the period-one unstable periodic orbit $\Gamma_+$, and Poincar\'e section $\Sigma$ defined as $x(t)-y(t) = 0$. In sub-figure (b) we plot the the projection of the $x$-component of the return map  of R\"ossler system. Assuming that a trajectory $(x(t),y(t),z(t))$ intersects with $\Sigma$ at $t = t_n$ for $n=1,2, ... $, we define $x_n = x(t_n) $. $(\gamma_x,\gamma_z)$ represents the intersection of the periodic orbit $\Gamma_+$ with the section $\Sigma$.}
   \label{fig:Ross_retmap}
\end{figure*}

\section{Weak tracking for nonautonomous R\"ossler system}
\label{sec:bif_Ross}

The R\"ossler system \cite{Rossler1976, Rossler1976a} proposed one of the simplest systems of ODEs that can have chaotic attractors. This has only one non-linear term and the system is given by:
\begin{equation}
\label{eq:RosslerSys}
\begin{array}{rcl}
\dot{x}&=&-y-z,\\
\dot{y}&=&x+ay,\\
\dot{z}&=&b+z(x-c). 
\end{array}
\end{equation}
There are many choices of parameters $a,b$ and $c$ that give chaotic attractors \cite{Alligood1996, Barrio2014, Cnrs1995}. We use as default $a=b=0.2$ and $c=5.7$ \cite{Rossler1976}, which give a chaotic attractor as shown in Figure~\ref{fig:Ross_retmap}(a).

We fix $b=0.2$ and $c=5.7$ throughout and analyse the bifurcations of \eqref{eq:RosslerSys} as $a$ varies between asymptotic values of $a_{\pm}$ as $t\rightarrow \pm\infty$. This ``frozen'' system has equilibria at 
$$
(x_{1,2},y_{1,2},z_{1,2})=\frac{c\pm\sqrt{c^2 - 4ab}}{2a} \left(a,-1,1 \right).
$$  
The equilibrium $p_1$ is asymptotically stable for any negative $a$ and bifurcates to stable periodic orbit at supercritical Hopf bifurcation point $a_{HB}\approx0.005978$. Soon after Hopf bifurcation, the resulting stable periodic orbit exhibits period doubling at $a_{\text{PD}}=0.1096$, and a period doubling cascade as $a$ increases until the system exhibit chaotic behaviour at $a\approx0.155$. 

To examine weak tracking, we shift $a$ from $a_{-}$ to $a_{+}$ for some $a_{-}, a_{+} \in \R$. Namely, 
$$
a(rt) =  \frac{\Delta}{2} \left( \tanh\left( \frac{\Delta r t}{2} \right)+1 \right) - a_{-}
$$ 
where $\Delta = a_{+} - a_{-}$, $r>0$ and $a_{-}$ ($a_{+}$) are the minimum (maximum) value of the parameter shift $a$. Throughout this paper we fix $a_{+}= -a_{-} = 0.2$. We can write the resulting R\"ossler system with parameter shift $a(t)$ as:  
\begin{equation}
\begin{array}{rcl}
\dot{x}&=&-y-z\\
\dot{y}&=&x+y \, a(rt)\\
\dot{z}&=&b+z(x-c) 
\end{array}
\label{eq:RosslerSysparshift}
\end{equation}

 The past limit system of \eqref{eq:RosslerSysparshift} has a hyperbolic stable equilibrium, $Z_- = \frac{c - \sqrt{c^2 - 4ba_-}}{2a_-} \left(a_-,-1,1 \right) \approx (-0.007, 0.0351, -0.0351)$. The future limit system, on the other hand, has a chaotic attractor $A_+$ that is the typical R\"ossler attractor in Figure~\ref{fig:Ross_retmap}(a). 
 
According to \cite[Theorem~2.2]{Ashwin2017},  for any $r>0$ system \eqref{eq:RosslerSysparshift} must have a pullback attracting solution $\A^{[a,r,Z_-]}$ that limits to $Z_-$, backward in time. Moreover, One can show that for almost every small enough $r>0$, the upper forward limit of the pullback attractor $\A^{[a,r,Z_-]}$ is the whole chaotic attractor $A_+$. Nevertheless, there is a set of isolated values of  $r>0$ that allow $\A^{[a,r,Z_-]}$ to end up tracking one of the unstable periodic orbits that are densely embedded in $A_+$. In this paper, we consider the period-one periodic orbit $\Gamma_+$, in particular, see Figure~\ref{fig:Ross_retmap}. However, similar arguments can be made for any unstable periodic orbits contained in $A_+$.

\begin{figure}  
\centerline{\includegraphics[width=0.8\columnwidth]{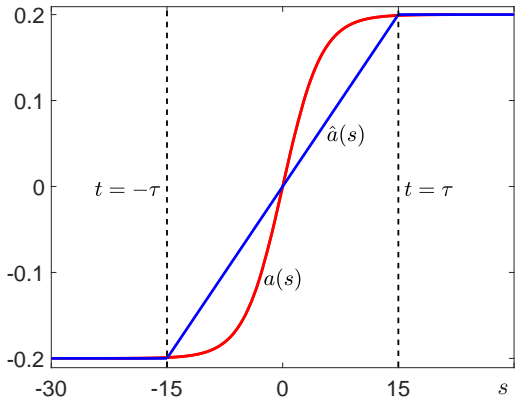}}
   \caption{The parameter shift $a(s)$ and the piecewise linear approximation $\hat{a}(s)$ vs time, for $a_{+} = - a_{-} =0.2$ and $\delta = 0.001$. }
   \label{fig:a_hat}
\end{figure}

\subsection{ Piecewise linear shift}
\label{sec:Piecewiseshift}

 In order to show that there are values of $r$ such that $\A^{[a,r,Z_-]}$ limits to $\Gamma_+$ as $t\to \infty$ we approximate the parameter shift $a(rt)$ by the following piecewise linear function $\hat{a}(rt)$:
\begin{gather*}
\hat{a}(s) = 
\begin{cases}
  a_{-}    &  ~~~~ s\in (-\infty , -\tau),\\
  (\Delta s + a_{+} + a_{-})/2 &  ~~~~ s\in [-\tau ,   \tau],\\
   a_{+}     &  ~~~~ s\in ( \tau  , \infty).
\end{cases}
\end{gather*}
where $\tau = \Big(\log({\Delta - \delta}) - \log(\delta)\Big)/\Delta$, for small enough $\delta>0$, note that at time $\pm \tau$ the value of $a$ is $\delta$-close to the upper and lower limits. i.e  $a(\tau)=a_{+} - \delta$ and $a(-\tau)=a_{-} + \delta$, see Figure~\ref{fig:a_hat}. 

The fact that $\hat{a}$ is fixed for any $t>\tau$, allows us to consider $A_+$ as an attractor for the system rather than just the upper forward limit of the pullback attractor $A_{t}^{[a,r,Z_-]}$.  

We embed a Poincar\'e section $\Sigma$ parametrised by $(x,z)$ with $x\leq 0$, as: 
$$
\{(x,x,z):~(x,z)\in\Sigma \}\subset \R^3,
$$

 and consider $t^*$, which is any real value that satisfies \begin{enumerate*} \item $t^* \geq \tau$ and \item $A_{t^*}^{[a,r,Z_-]} \in \Sigma$ \end{enumerate*}, i.e. $A^{[\hat{a},r,Z_-]}_{t^*}$ is a point in $\Sigma$.

 Note that, the intersection of $\Gamma_+$ with $\Sigma$ is a fixed point $\gamma$ for the return map. If $r_c>0$ is chosen such that, $A^{[\hat{a},r,Z_-]}_{t^*}$ is one of the pre-images of $\gamma$, then the the upper forward limit of $\A^{[\hat{a},r,Z_-]}$ is $\Gamma_+$ and $r_c$ is a critical rate for weak tracking.

\begin{figure*}
\centerline{\includegraphics[width=1.8\columnwidth]{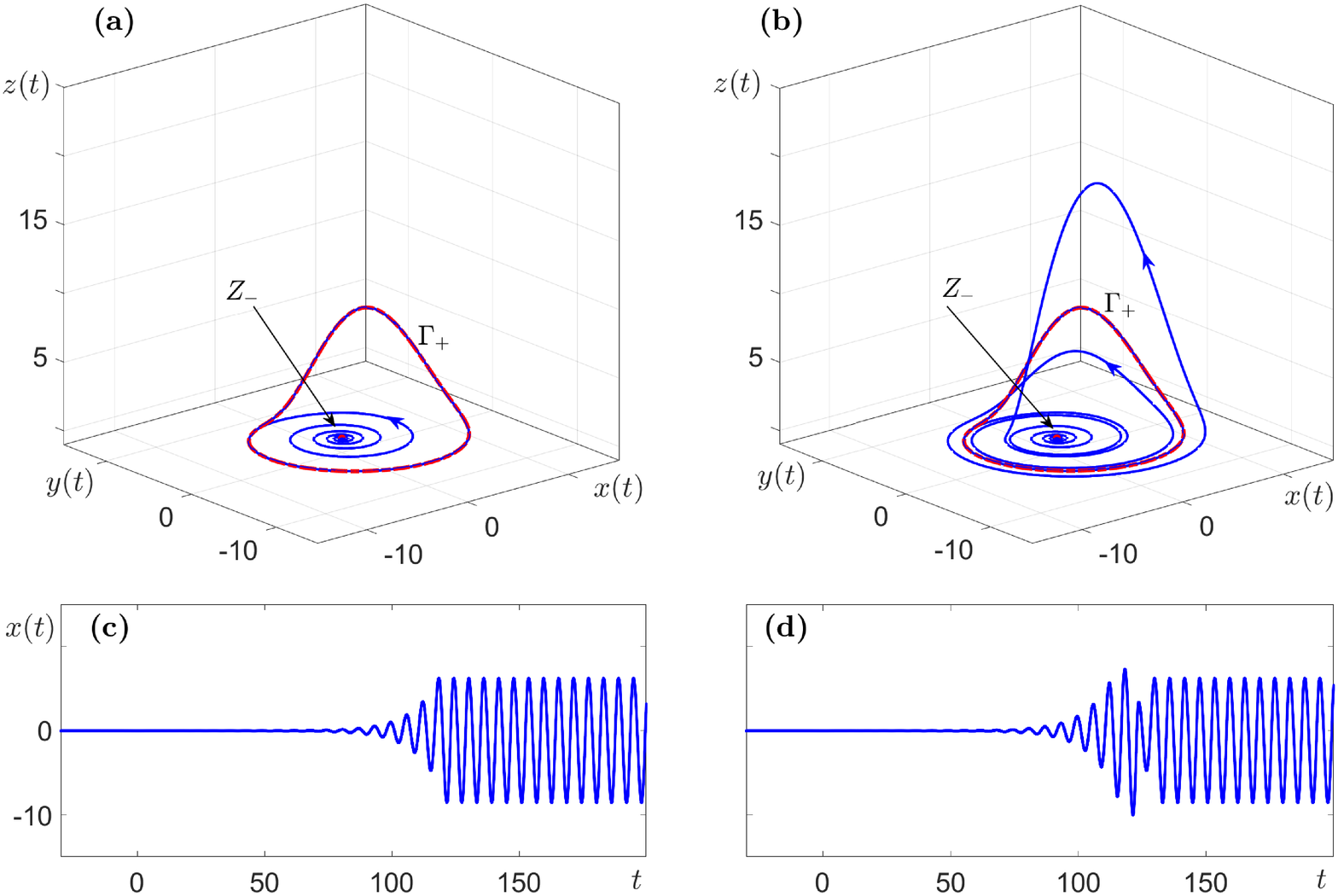}}   
 \caption[Two examples of weak tracking (EtoP connection) for \eqref{eq:RosslerSysparshift}]{Two examples of weak tracking (EtoP connection) for \eqref{eq:RosslerSysparshift}. The parameters are $b=a_{\max}=-a_{\min}=0.2$, $c=5.7$ and $T=150$. (a) and (c) show the EtoP connection at $r=0.9202212159423$, (b) and (d) show the connection at $r=0.995651959127$. } 
 \label{fig:Ross_EtoPconn}
\end{figure*}

\subsection{ Density of critical rates: Numerical evidence}
\label{sec:numarics}

\begin{figure}
\centerline{\includegraphics[width=1\columnwidth]{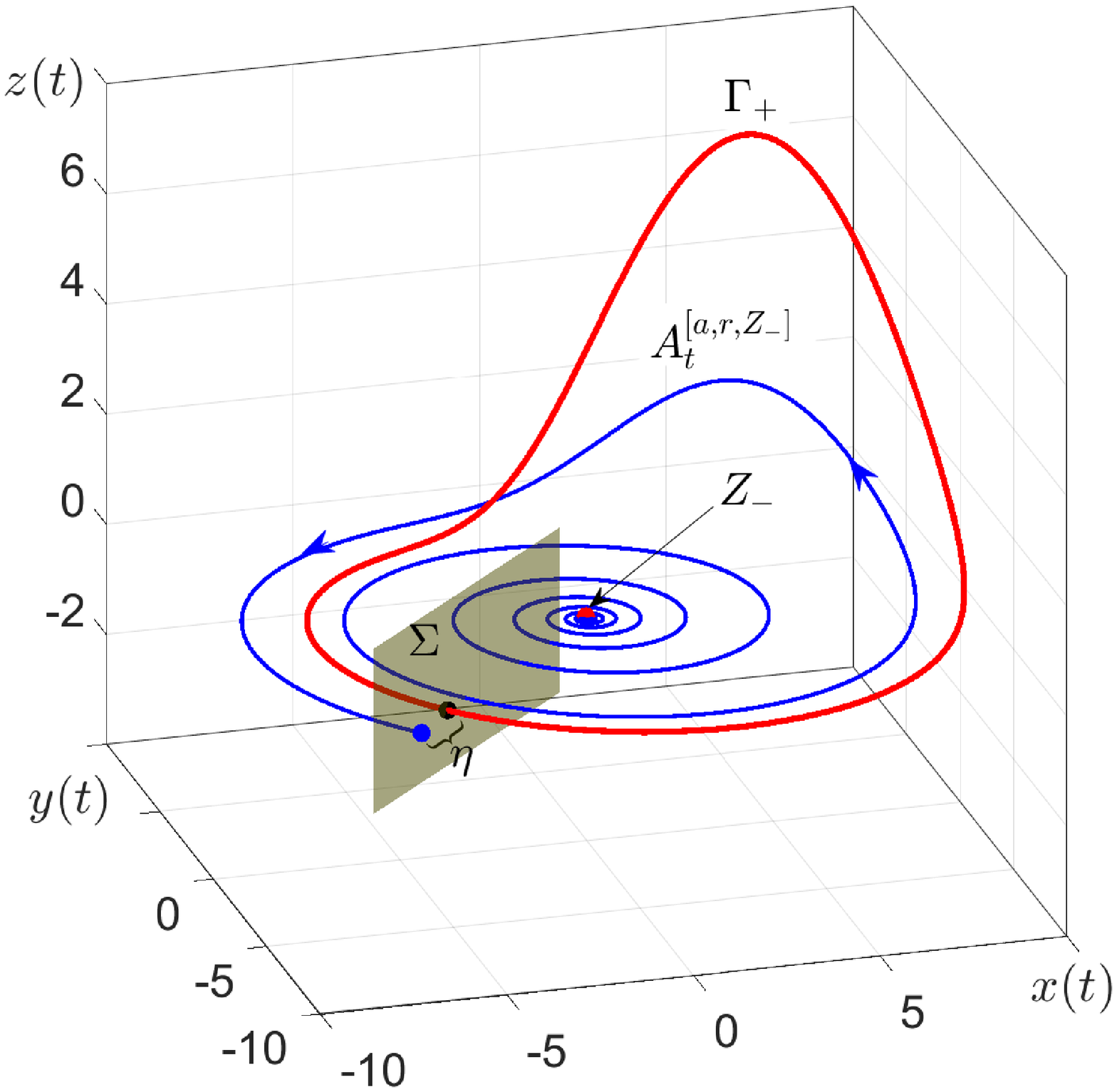}}
\caption{ A schematic diagram showing the shooting method we use to find the the connection between $Z_-$ and $\Gamma_+$ for \eqref{eq:RosslerSysparshift}, see Appendix~\ref{sec:supp} for animated version of this figure.} 
\label{fig:Ross_schdiag} 
 \end{figure}

To investigate weak tracking for System~\eqref{eq:RosslerSysparshift}, with the smooth parameter shift $a(rt)$, we use a shooting method as follows:
\begin{enumerate}

\item We approximate the pullback attractor $\A^{[a,r,Z_-]}$ by integrating \eqref{eq:RosslerSysparshift}, subject to an initial condition $Z_{\text{init}}$ fairly close to $Z_-$. Namely, we choose $Z_{\text{init}} = (-0.007,0.035,-0.035)$ and the integration time is from $-30$ to $T$.
\item The point pullback attractor can be given as $\A^{[a,r,Z_-]} = (\tilde{x}^{r}(t),\tilde{y}^{r}(t),\tilde{z}^{r}(t))$, where $t\in[-30,T]$. 
\item Recall that the Poincar\'e section $\Sigma$ is parametrised by $(x,z)$ with $x\leq 0$, as: 
$$
\{(x,x,z):~(x,z)\in\Sigma \}\subset \R^3
$$
\item Assume that $\A^{[a,r,Z_-]}$ intersects $\Sigma$ at times $t_n\leq T$ for $n = 1,2, ... ,N$, $N\in\mathbb{N}$ and $t_{n-1}<t_{n}$. 
\item Consider the final intersection point  $(\tilde{x}^r(t_N),\tilde{z}^r(t_N))\in\Sigma$. We approximate a signed distance from the stable manifold of $\gamma$ by the following real valued ``gap function''
 $$
 \eta(r):= \frac{\big((\tilde{x}^r(t_N),\tilde{z}^r(t_N)) - \gamma \big) v_{\text{s}}^T}{v_{\text{s}}v_{\text{s}}^T},
 $$ 
where $\gamma=(\gamma_x,\gamma_z) \in \Sigma$ is the fixed point of R\"ossler return map, see Figure~\ref{fig:Ross_retmap}, and $v_{\text{s}}$ is stable eigenvector of $\gamma$ for the return map. Note that $\eta(r)$ also depends on  $T$, $b$, $c$, $a_{\min}$, $a_{\max}$ and $Z_{\text{init}}$. However, here we only consider variation of $r$.  
\item By analogy to Section~\ref{sec:Piecewiseshift}, whenever $\eta(r_c) \approx 0$ the pullback attractor $\A^{[a,r,Z_-]}$  intersects the stable manifold of $\Gamma_+$, which gives the desired EtoP connection. In other words, $\A^{[a,r,Z_-]}$ weakly tracks $A_+$ at $r=r_c$. The method is illustrated schematically in Figure~\ref{fig:Ross_schdiag}.
\end{enumerate}

The function $\eta(r)$ is as smooth as the state variables of \eqref{eq:RosslerSys}, i.e. it is at least $C^1$. Consequently, one can numerically approximate its roots, and hence the critical rates of weak tracking, using a root-finding algorithm such as Newton-Raphson method. Figure~\ref{fig:Ross_EtoPconn} shows that system \eqref{eq:RosslerSys} exhibit weak tracking at two different critical rates.

\begin{figure} 
\centerline{\includegraphics[width=0.9\columnwidth]{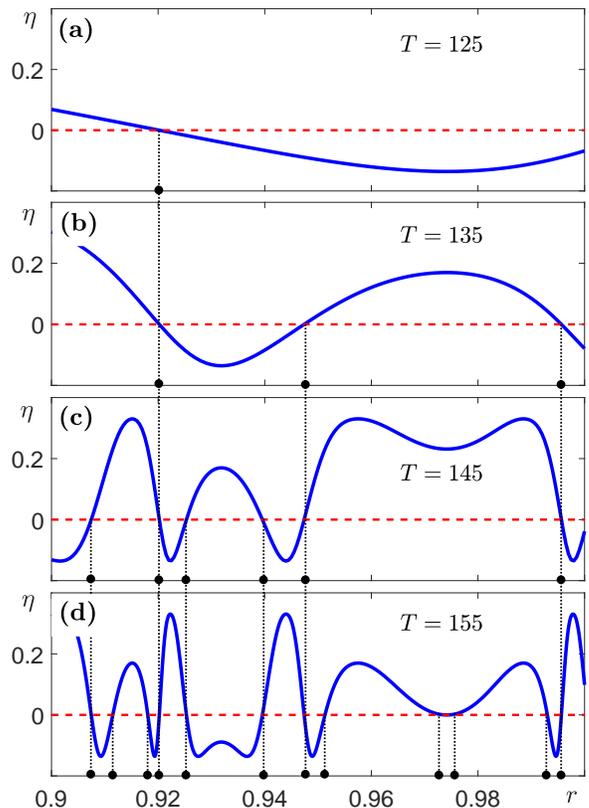}}
    
   \caption{Graphs of $\eta$ for increasing integration time $T$ to add additional intersections of $\Sigma$: roots correspond to connections from $A_-$ to the periodic orbit $\Gamma_+$. (a) $T=125$, (b) $T=135$, (c) $T=145$, and (d) $T=155$. It can be seen that additional zeros of $\eta$ (corresponding to critical rates that give weak tracking) appear as $T$ increases. The parameter values are $b=a_{\max}=-a_{\min}=0.2$ and $c=5.7$. }
\label{fig:Ross_etaroots}
\end{figure}

We point out two numerical difficulties in our numerical approach to approximate the rates of weak tracking: First, there is a large delay in Hopf bifurcation that forces us to choose fairly large integration time $T$ in our calculations, which increases the computational cost. Delay in dynamic bifurcations is very common and not easy to avoid. For a system with linearly changing time-dependent parameter with slope $r$, dynamic Hopf bifurcation may have a delay time proportional to $1/r$ before fast escape from the curve of unstable equilibria occurs \cite{Neishtadt2009, Lobry1991}. More details on dynamic bifurcations and their delay can be found in \cite[Chapter~2]{Berglund2006}. Second, Figure~\ref{fig:Ross_etaroots} shows that $\eta(r)$ is smooth with respect to $r$ for a particular range of $r$, which is $[0.9,1]$. However, there is no guarantee that $\eta(r)$ is smooth or even continuous for finite $T$. The definition of $\eta(r)$ depends on the maximum intersection time which in turn depends on the integration time $T$. Nevertheless, $T$ can be chosen to smooth $\eta(r)$ out for any range of $r$.

Our numerical investigation suggests that there are infinitely many critical rates that give weak tracking for \eqref{eq:RosslerSys}. In Figure~\ref{fig:Ross_etaroots} we plotted  $\eta(r)$ against $0.9 \leq r\leq 1$, for different values of $T=125$, $135$, $145$ and $155$. The results show that as $T$ increases, the number of roots of $\eta(r)$ increases rapidly.  

Despite the other periodic orbits that are embedded in in $A_+$, even for just one periodic orbit $\Gamma_+$, our numerical investigation shows that there are infinitely many critical rates that give weak tracking. In fact, we believe that the set of all critical rates $r_c$ is dense in $\R$.

\section{Discussion}
\label{sec:discussion}

We study the well known R\"ossler system \eqref{eq:RosslerSysparshift} with parameter shift, as a tool to illustrate a new rate-induced phenomenon that we term ``weak tracking''. We monotonically shift the bifurcation parameter $a$ such that the system has an equilibrium attractor for the past limit system and chaotic attractor for the future limit system. We then show that there are isolated critical rates at each of which the pullback attractor solution of the system ends up tracking an embedded saddle periodic orbit in the future chaotic attractor. We use a numerical approach, based on shooting method and carefully chosen Poincar\'e section, to approximate these critical rates.

For the nonautonomous  R\"ossler system \eqref{eq:RosslerSysparshift} with a parameter shift from stable equilibrium to chaos, we suggest there is a dense set of critical rates that give weak tracking. We give an argument below that this is the case if the system has piecewise linear forcing instead of smooth parameter shift and provide in Figure~\ref{fig:Ross_etaroots} numerical evidence of the existence of the dense set of critical rates for smooth parameter shift. 

Although our example considers a specific choice of parameters, the necessary ingredients for weak tracking are present in a wide range of the parameter space of the nonautonomous R\"ossler system. These ingredients are simply (i) a hyperbolic attracting equilibrium for the past limit system (ii) a chaotic (non-minimally invariant) attractor for the future limit system and (iii) a rate dependent shift in parameters that means for certain rates the pullback attractor gets ``caught'' in unstable dynamics within the chaos.

More precisely, in order for the parameter shift system \eqref{eq:parshift} to exhibit weak tracking along a branch of attractors $A(\lambda)$ from a past limit attractor $A_-$ to $A_+$, it is clear that the future limit system must have a proper invariant subset $S_+$ (in our case we consider $S_+=\Gamma_+$) of the future limit attractor $A_+$, and the pullback attractor with past limit $A_-$ must ``fit in" to $S_+$. If we consider (\ref{eq:parshift}) then weak tracking corresponds to existence of a pullback attractor $A^{[\Lambda,r,A_-]}_t$ with backward limit $A_-$ and forward limit $S_+$. This will only be possible if the dimension of $A_-$ is small enough with respect to that of $S_+$. For an eventually constant parameter shift such as in Fig.~\ref{fig:a_hat}, note that $A^{[\Lambda,r,A_-]}_t=A_-$ as long as $t$ is sufficiently negative, and as nonautonomous time evolution will be a diffeomorphism between any two finite times, i.e. $A^{[\Lambda,r,A_-]}_t$ is diffeomorphic to $A_-$ for all finite $t$. Hence in this eventually constant case a necessary condition for $A^{[\Lambda,r,A_-]}_t$ to limit to $S_+$ is that $A^{[\Lambda,r,A_-]}_t\subset W^s(S_+)$ for sufficiently large $t$ where $W^s$ is the stable set for the future limit flow. Hence 
$$
\dim(A_-)=\dim\left(A^{[\Lambda,r,A_-]}_t\right)\leq \dim(W^s(S_+))
$$
(where $\dim(A)$ represents Hausdorff dimension of $A$). Hence weak tracking require
\begin{equation}
\dim(A_-)\leq \dim(W^s(S_+))
\end{equation}
which means in particular if $\dim(A_-)>\dim(S_+)$ then a connection is not possible.

Moreover, note that for large enough $t$, the set $A^{[\Lambda,r,A_-]}_t$ will, in the generic case, vary nontrivially with $r$. Any interaction between this and $W^s(S_+)$ will typically be transverse on varying $r$: this argues that values of $r$ where there is weak tracking are isolates. Density of $W^s(S_+)$ within the basin $\mathcal{B}(A_+)$ of $A_+$, with respect to the future limit flow, would imply the density of a set of critical rates giving weak tracking to this $S_+$.

For example, if $A_-$ is an equilibrium or periodic orbit then it is possible to have weak tracking to a periodic orbit $S_+$ contained in $A_+$ a chaotic attractor. If $A_-$ is chaotic then weak tracking will only be possible to an invariant set $S_+$ with dimension greater than $A_-$. A similar result will presumably apply more generally, even if the shift is not eventually constant. In this case the condition for weak tracking will be in terms of a condition for existence of a connection from $A_-$ to $S_+$ for the extended autonomous system.

Parameter shift systems such as (\ref{eq:parshift}), and asymptotic autonomous systems more generally, have a rich tipping behavior. Ref. \cite{Alkhayuon2018} gives an example of a system with pullback attractor that exhibit partial rate-dependent tipping, where an entire subset of the pullback attractor tracks different quasi-static attractor than it would be for other rates of shift, while the rest of the pullback attractor still tracks the associated quasi-static attractor. This behaviour can still be produced in R\"ossler system with a suitable parameter shift that shifts the chaotic attractor partially out of its basin of attraction. 

More precisely, suppose we have a parameter shift $\Lambda(rt)$ that limits to $\lambda_{\pm}$ forward and backward in time respectively, such that the attractors for the future and the past limit systems, $A_{\pm}$, are non-equilibrium attractors. Ref.  \cite{Alkhayuon2020} shows that partial tipping is possible, for some values of $r$, if:
\begin{equation}
\begin{aligned}
&A_- \not\subset \mathcal{B}(A_+), 
\end{aligned}
\end{equation} 

Besides the phenomena illustrated in Ref. \cite{Alkhayuon2018}, nonautonomous systems with nonequilibrium attractors may exhibit other transitions. For example, systems that have attractors with fractal basin boundaries may exhibit fractality-induced tipping \cite{Balint} due to the high complexity of the basin not because of the well known tipping mechanisms presented in Ref. \cite{Ashwin2011}. Basins of attraction with fractal boundary are very common in physical systems, and can cause a high uncertainty when it comes to predicting the final state of a trajectory. We refer to Ref. \cite{Aguirre2009} for further details. Fractal boundaries can result from crossings of the stable and unstable manifold of an invariant set that is embedded in basin boundary. 

Fractality may also be a sign of the presence of {\em transient chaos} \cite{Tel2015}. One phenomenon that can lead to transient chaos is a {\em boundary crisis} \cite{Grebogi1983a, Grebogi1982}, where the attractor intersects its basin boundary and leaks out. If the time dependent parameter passes through a region where there is a crisis, then the system exhibit {\em attractor hopping} behaviour \cite{kraut2002}, which may led to partial or even total tipping. 
 
~
 
\section*{acknowledgements}
HA's research is funded by Enterprise Ireland and Laya Healthcare, grant agreement No. 20190771. PA’s research was partially supported by the CRITICS Innovative Training Network, funded by the European Union’s Horizon 2020 research and innovation programme under the Marie Sklodowska-Curie Grant Agreement No. 643073. Both authors gratefully acknowledge the insightful comments by Jan Sieber, Roberto Barrio, Damian Smug, Paul Richie, and Sebastian Wieczorek.

~

\appendix
\section{Supplementary material}
\label{sec:supp}
We provide a MATLAB code for the shooting method we used in Section~\ref{sec:numarics} and a supplementary movie shows animations of Figure~\ref{fig:Ross_etaroots} in the GitHab repository: ``github.com/hassanalkhayuon/WeakTracking''.
\bibliography{../bibliography}

\end{document}